\documentclass[letterpaper, 10 pt, conference]{ieeeconf}  % Comment this line out if you need a4paper

\IEEEoverridecommandlockouts                              % This command is only needed if
\usepackage{epsfig} % for postscript graphics files
\usepackage{mathptmx} % assumes new font selection scheme installed
\usepackage{times} % assumes new font selection scheme installed
\usepackage{amsmath} % assumes amsmath package installed
\usepackage{amssymb}  % assumes amsmath package installed
\usepackage{epstopdf}
\usepackage{color}
\usepackage{cite}
\usepackage{subfigure}

\newtheorem{lemma}{Lemma}

\newtheorem{assumption}{Assumption}
\newtheorem{proposition}{Proposition}

\title{\LARGE \bf
Performance Regulation and Tracking via Lookahead Simulation: Preliminary Results and Validation
}

\author{Y. Wardi$^{*}$, C. Seatzu$^{**}$, M. Egerstedt$^{*}$, and I. Buckley$^{*}$
\thanks{$^*$School of Electrical and Computer Engineering, Georgia Institute of Technology, Atlanta, GA 30332. Email:
ywardi@ece.gatech.edu, magnus@ece.gatech.edu, ihbuckl@g.clemson.edu.}
\thanks{$^{**}$Department of Electrical and Electronic
Engineering, University of Cagliari,  Italy. Email: seatzu@diee.unica.it.}
}

\begin{document}

\maketitle
\thispagestyle{empty}
\pagestyle{empty}

\begin{abstract}

  This paper presents an approach to target tracking that is based on a variable-gain integrator and
the Newton-Raphson method for finding  zeros  of a function. Its underscoring idea is the determination of the feedback law
by measurements of the system's output and estimation of its future state via lookahead simulation.
The resulting feedback law is generally nonlinear. We first apply the
proposed  approach to tracking a constant reference by the output of nonlinear memoryless plants. Then
 we extend it in a number of directions, including the tracking of time-varying reference signals by
 dynamic, possibly unstable systems. The approach is new  hence its  analysis is preliminary, and theoretical results
 are derived  for nonlinear memoryless plants and linear dynamic plants.
 However, the setting for the controller does not require  the plant-system to be either linear or stable, and this is verified
 by simulation of an inverted pendulum tracking a time-varying signal.
  We also demonstrate results of laboratory experiments of controlling  a platoon of mobile robots.
 \end{abstract}

\section{Introduction}
Integrative action is an essential element in the steady-state tracking of a constant reference signal by the output of a linear
system. However, it is well known that a controller comprised solely of an integrator may have destabilizing
effects on
the closed-loop system. Therefore tracking controllers often include  proportional and derivative elements in addition to an
integrator \cite{Franklin14}, as is standard fare in any undergraduate controls classes. However, purely integrative actions have been used, for example for the regulation of computer processors; Ref.  \cite{Almoosa12} proposed
a standalone integral controller endowed with a variable gain in order to enhance the stability margins of the closed-loop system.\footnote{These applications
include the regulation of power and instruction throughput by the processor's clock frequency. The rationale behind the choice
of such a control architecture, as well as results of simulations  on industry benchmark programs, and implementations on Intel's Haswell microarchitecture
\cite{Hammarlund14},  can be found in \cite{Almoosa12, Chen15, Chen16}.}

This idea of a variable gain integrator can be generalized. To illustrate this idea, consider, for example, the
single-input-single-output discrete-time system in Figure 1. The objective of the controller is to have the system's output, $y_n\in R$,  with $n=1,2,\ldots$, denoting time,
asymptotically track the given reference  $r\in R$. Let $e_n:=r-y_n$ denote  the error signal, and let $u_n\in R$ be the input signal
(control signal) to the plant. The variable gain integrator takes on the form
\begin{equation}
u_n=u_{n-1}+A_n e_{n-1},
\end{equation}
where $A_{n}>0$ is its gain at time $n$. Note that if $A_{n}=A$ $\forall n=1,\ldots$, for a given constant
$A$, then the controller acts as an adder, the discrete-time equivalent of an integrator, and hence we call it
a variable-gain integrator.  Generally $A_{n}$ is not a constant.

Suppose that the plant is described by a memoryless nonlinearity,
$y_n=g(u_n)$ for a continuously-differentiable function $g:R\rightarrow R$. Then  $A_n$ can be   defined as follows,
\begin{equation}
A_{n}:=\Big(\frac{\partial g}{\partial u}(u_{n-1})\Big)^{-1}\big(r-y_{n-1}\big),
\end{equation}
where we assume that the derivative term in Eq. (2) exists and is non-zero.

The resulting tracking (regulation) algorithm consists of a recursive application of Eq. (2), and we note that it comprises an implementation of
the Newton-Raphson method for solving the algebraic equation $g(u)=r$.
In \cite{Almoosa12,Chen15,Chen16} the use of this variable gain integrator was argued to be competitive
with extant techniques for controlling power and instruction-throughput in multicore computer processors despite
its simple form.
An analysis in a general setting of nonlinear, memoryless systems was carried out in \cite{Wardi16}.

\vspace{.2in}
\begin{figure}[h]
\centering
\vspace{-.3cm}
\includegraphics[width=0.40\textwidth]{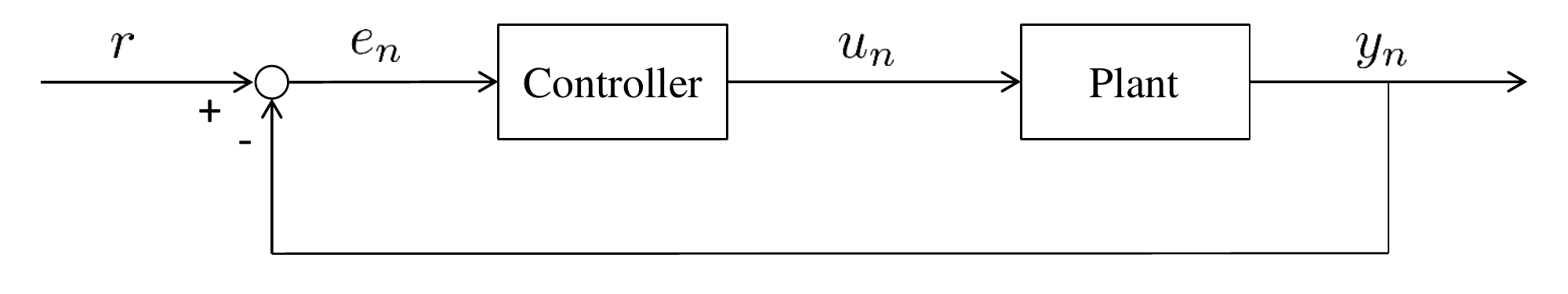}
% %\includegraphics[width=3.15in,angle=0]{Adaptive.jpg}
{\small \caption{Basic Control System}}
\end{figure}

The objective of this paper is to extend the aforementioned regulation  approach from discrete-time memoryless systems to continuous-time dynamical systems,  and from the regulation of constant reference values to
the tracking of time-varying reference signals.
The extension from discrete-time systems to memoryless continuous-time
systems is  straightforward
and consists of replacing Eq. (2) by a differential equation. However, an  extension of the plant from a memoryless system to a dynamical system is more subtle. We propose such an extension by using a nonlinear observer,      based on
a lookahead simulation of the plant. A key question is how to choose the time-horizon
for the lookahead simulation. Large time horizons  may yield tracking convergence only
for constant reference signals, while the tracking of time-varying signals requires short time horizons. However,
short time horizons may render the closed-loop system unstable.
 To get around this problem we speed up the controller subsystem  (but not the plant)  which can  restore stability  thereby yielding tracking convergence.

To our knowledge this approach to tracking is new, and the main objective of the paper is to introduce and explain it, derive preliminary theoretical
results, and present results of simulation and laboratory experiments. The general theoretical problem is to classify the systems for which the
proposed controller is provably convergent. In such generality the problem is beyond the scope of this paper, but we derive
convergence results for memoryless nonlinear plants and dynamic linear plants.
Furthermore, we  analyze in detail the particular example
of position control in systems obeying Newton's second law with a drag. The presented simulation results run ahead of the theoretical developments,
and they  verify the principles of the tracking technique for
various linear and nonlinear systems, with both stable and unstable plants, and with  constant target levels as
well as time-varying reference signals.

It should be pointed out that the issue of nonlinear regulation is certainly a well-established topic and in this paper we do not provide solutions to problems that were previously not solved. In particular, techniques such as the Byrnes-Isidori regulator \cite{Isidori90} based on internal model techniques, or Khalil's high-gain observers for output regulation \cite{Khalil98}, certainly are more powerful than the methods introduced in this paper.
However, the effectiveness of these regulators rely on significant computational sophistication, such as nonlinear inversions and the appropriate nonlinear normal form, e.g., \cite{Isidori95,Sastry99}. As such, the contribution in this paper should be understood as a computationally unproblematic, variable gain integrator, and an initial exploration as to when such a structurally simple controller can indeed achieve the desired performance.

In the  rest of the paper  Section 2 presents the regulation technique with its observer and derives
theoretical results, Section 3 provides simulation and experimental verifications, and Section IV concludes the paper.

\section{Regulation and Tracking}

We start this section by extending the discussion of the system depicted in Figure 1 from
discrete time to continuous time.  Consider the feedback system shown in Figure 1 except that    time $t\geq 0$ is continuous, and hence
the signals around the loop are $u(t),~y(t)$, and $e(t)$. Suppose that all of these signals  are $k$-dimensional
for some $k\geq 1$,  and the reference $r$ is a $k$-dimensional function of time,  $r(t)\in R^k$.

\subsection{Memoryless plant}
Consider first the case where the plant is
represented by a function $g:R^k\rightarrow R^k$, such that for every $t\geq 0$,
\begin{equation}
y(t)=g(u(t)).
\end{equation}
Suppose also that the function $g(u)$ is continuously differentiable, and that its Jacobian
$\frac{\partial g}{\partial u}(u)$ is nonsingular for every $u\in R^k$ considered in the sequel.

A natural continuous-time equivalent equation to Eqs. (1) and (2) is
\begin{equation}
\dot{u}(t)=\Big(\frac{\partial g}{\partial u}\big(u(t)\big)\Big)^{-1}\big(r(t)-g(u(t)\big);
\end{equation}
here $u(t)$ moves in the direction defined by Newton-Raphson algorithm but the step size is scaled by $dt$.
\begin{proposition}
The following limit is in force:
\begin{equation}
\limsup_{t\rightarrow\infty}||r(t)-y(t)||~\leq~\sup\big\{||\dot{r}(t)||~:~t\geq 0\big\}.
\end{equation}
\begin{proof}
Define $\rho:=\sup\big\{||\dot{r}(t)||~:~t\geq 0\big\}$. Define the function
$V:R^k\times R^+\rightarrow R^+$ by
\begin{equation}
V(u,t)=\frac{1}{2}||r(t)-y(t)||^2.
\end{equation}
Taking derivatives with respect to $t$, and by Eqs. (3) and (4),
\begin{equation}
\dot{V}(u(t),t)=\big(r(t)-y(t)\big)^{\top}\big(\dot{r}(t)-(r(t)-y(t))\big).
\end{equation}
Fix $\varepsilon\in(0,1)$. By Eq. (7), if
$||r(t)-y(t)||>(1+\varepsilon)\rho$, then
\begin{equation}
\dot{V}(u(t),t)~\leq~-\varepsilon(1+\varepsilon)\rho^2.
\end{equation}
By Lyapunov direct method \cite{Khalil02} it follows that $\limsup_{t\rightarrow\infty}||r(t)-y(t)||\leq(1+\varepsilon)\rho$,
and since $\varepsilon>0$ can be arbitrarily small, Eq. (5) is satisfied.
\end{proof}
\end{proposition}

{\it Remarks:} 1). In the special case where $r(t)=r\in R^k$, a constant, the function $V(u):=V(u,t)$ is a Lyapunov function,
and it follows that, as $t\rightarrow 0$, $\lim y(t)=r$.

2). Proposition 1 means that asymptotically $y(t)$ is confined to a ball with center $r(t)$ and radius $\rho$. To reduce
this radius we can increase the controller's gain by multiplying the RHS of (4) by $\alpha>1$. The control equation
then becomes
\begin{equation}
\dot{u}(t)=\alpha\Big(\frac{\partial g}{\partial u}\big(u(t)\big)\Big)^{-1}\big(r(t)-g(u(t)\big).
\end{equation}
With the  Lyapunov function $V(u,t)$ defined by (6), the same arguments comprising  the proof of Proposition 1  yield that
\begin{equation}
\limsup_{t\rightarrow\infty}||r(t)-y(t)||~\leq~\frac{1}{\alpha}\sup\big\{||\dot{r}(t)||~:~t\geq 0\big\}.
\end{equation}
We next extend this tracking control law to dynamic plants.\\

\subsection{The plant as a dynamical system}
Suppose that the plant is represented by a differential equation of the form
\begin{equation}
\dot{x}(t)=f(x(t),u(t)),~~~~~~~t\geq 0,
\end{equation}
where $x(t)\in R^n$,  for a function $f:R^n\times R^k\rightarrow R^n$;
the
 initial state is  $x_0:=x(0)$. The output $y(t)$  is given by
\begin{equation}
y(t)=h(x(t))
\end{equation}
for a function $h:R^n\rightarrow R^k$.

The following assumption ensures that the state variable $x(t)$  and the output $y(t)$ are well defined
for all $t\geq 0$.
\begin{assumption}
1). For every $u\in R^k$, the function $f(x,u)$ is  continuously-differentiable in
$x$, and the functions $f(x,u)$ and  $\frac{\partial f}{\partial u}(x,u)$  are locally Lipschitz continuous in $(x,u)$.
2). The function $h(x)$ is continuously differentiable in $x$. 3). For every compact set
$U\subset R^k$ there  exist $K_1>0$ and $K_2>0$ such that, for every $(x,u)\in R^k\times U$,
$||f(x,u)||\leq K_1|x||+K_2$.
\end{assumption}

In extending the control law from the memoryless case (Eq. (4)) to the dynamic-plant case it is necessary  to first
define the function
$g(u)$, and there is no single natural way to do it. Our choice of $g(u)$ is based on an evaluation of the system's output
$T$ seconds in the future by a lookahead simulation for a given $T>0$.  Specifically, suppose that at time $t$ the
state $x(t)$ is measured, and a lookahead simulator computes what the output would be at time $t+T$ if the input
were  $u(\tau)=u(t)$ $\forall\tau\in[t,t+T]$. The result of this computation is $g(u)$. Formally,
given $t\geq 0$, $x\in R^n$,   $u\in R^k$ and $T>0$, let $\xi(\tau)\in R^n$ be defined by the equation
\begin{equation}
\dot{\xi}(\tau)=f(\xi(\tau),u)
\end{equation}
with the boundary condition $\xi(t)=x$; note that the input is a constant $u(\tau)=u$ $\forall\tau\in[t,t+T]$.
We use the notation $\phi(x,t;u,\tau)=\xi(\tau)$ to emphasize its dependence on $x,~t$ and $\tau\geq t$ in addition to $u$.
In a specific run of the system, let $x(t)$ denote its state variable at time $t$ so that its output is $y(t)=h(x(t))$ according to Eq. (12). Suppose also
that at time $t$, $x(t)$ can be measured and a lookahead simulation can compute $\phi(x(t),t;u(t),t+T)$; note that in this simulation the
input to the plant-model in the simulator is $u(t)$ $\forall\tau\in[t,t+T]$. Then $g(u)$ is defined as follows,
\begin{equation}
g(u)=h(\phi(x(t),t;u(t),t+T)).
\end{equation}

Of course  $g(u)$ is a function of $t$, $x(t)$ and
$T$ in addition of $u$, but we use the simplified notation $g(u)$ when no confusion
arises.

The feedback law is defined by adjusting Eq. (4) to the present case of dynamic plants. The
only difference is to replace $r(t)$ in (4) by $r(t+T)$, resulting in the equation
\begin{equation}
\dot{u}(t)= \Big(\frac{\partial g}{\partial u}(u(t))\Big)^{-1}\big(r(t+T)-g(u(t))\big).
\end{equation}
The reason for this change is that Eq. (4) is designed to match $g(u(t))$ to $r(t)$, whereas in Eq. (15) $g(u)$ is an estimator
of $y(t+T)$ and hence it is designed to match $r(t+T)$.

The closed-loop system is comprised of the plant and the controller subsystems, defined by Eqs. (11) and
 (15), respectively.  The computation of (15) including its various ingredients,  $g(u(t))$ and $\frac{\partial g}{\partial u}(u(t))$,
 can be performed by any numerical technique for differential equations; we use the forward-Euler method in all examples
 described in the sequel.

 Thus far the discussion of Eq. (15) has been predicated on its exact computation. However, in practical situations one can expect errors
 due to several factors, including measurements of $x(t)$, modelling uncertainties, and computational errors. An
 error analysis cannot fit in the paper due to its space limitation and hence will be presented in a forthcoming publication. We anticipate
 a result
 similar to Eq. (5) where $\dot{r}(t)$ is replaced by a cumulative measure of the various errors.

An extension of Proposition 1 to the present setting of dynamical systems does not work for the following two reasons: First, the closed-loop
system may be unstable, while this is not a problem for memoryless plants as can be seen in the proof of Proposition 1.
Second, such a result, if true, would imply the tracking of $r(t+T)$ by $g(u)$ which is not the same as
$y(t+T)$. Generally the error term $||g(u)-y(t+T)||$ can be made small by choosing a small $T>0$. However, we shall see that
often the closed-loop system is stable for large $T$ and unstable for small $T$. One approach to this difficulty is to choose
a small $T$ and then try to stabilize the system by scaling up  the Right-Hand Side (RHS) of Eq. (15). This is not
the same as increasing the controller's gain (unless the plant system is memoryless)  and
may have a stabilizing effect on the closed-loop system.   We shall see that for a particular class of Linear, Time-Invariant  (LTI) systems,
for a given value of $T>0$, the
tracking   of a constant reference
  $r\in R^k$ is achieved as long as the closed-loop system is stable.

The derivation of sufficient conditions for tracking for general nonlinear systems is beyond the scope of this
paper. Instead, we  next derive some theoretical results only for linear systems, but later provide simulation results
for linear and nonlinear systems.

\subsection{The plant as an LTI system}

Consider the case where the plant has the form
\begin{equation}
\dot{x}(t)=Ax(t)+Bu(t),~~~~~~y(t)=Cx(t),
\end{equation}
where the matrix-dimensions are $A\in R^{n\times n},~ B\in R^{n\times k},~ C\in R^{k\times n}$.
Suppose that $A$ is nonsingular. Fix $T>0$, and suppose that the matrix
$\big(CA^{-1}(e^{AT}-I_n)B\big)$ is nonsingular as well, where $I_n$ denotes the $n\times n$ identity
matrix.
By solving the differentiable equation (13), and by Eq. (14), it is readily
seen that
\begin{equation}
g(u)=C\Big(e^{AT}x(t)+A^{-1}(e^{AT}-I_n)Bu\Big),
\end{equation}
and hence
\begin{equation}
\frac{\partial g}{\partial u}(u)=CA^{-1}(e^{AT}-I_n)B.
\end{equation}
By Eq. (15), after some algebra it follows that
\begin{equation}
\dot{u}(t)=\Big(CA^{-1}(e^{AT}-I_n)B\Big)^{-1}\big(r(t+T)-Ce^{AT}x(t)\big)-u.
\end{equation}
Define the $(n+k)\times(n+k)$ matrix
$\Phi_T$ by
\begin{equation}
\Phi_{T}=\left(
\begin{array}{cc}
A & B\\
-\Big(CA^{-1}(e^{AT}-I_n)B\Big)^{-1}Ce^{AT} & -I_k
\end{array}
\right),
\end{equation}
where $I_k$ is the $k\times k$ identity matrix,
and define the $k\times k$ matrix $\Psi_T$ by
\begin{equation}
\Psi_T=\Big(CA^{-1}(e^{AT}-I_n)B\Big)^{-1}.
\end{equation}
Then
 the closed-loop system comprised of Eqs. (16) and (19) has the form
\begin{equation}
\left(
\begin{array}{c}
\dot{x}(t)\\
\dot{u}(t)
\end{array}
\right)
=
\Phi_T
\left(
\begin{array}{c}
x(t)\\
u(t)
\end{array}
\right)
+
\left(
\begin{array}{c}
0\\
\Psi_T
\end{array}
\right)r(t+T).
\end{equation}

\begin{lemma}
Suppose that $A$ is Hurwitz. There exists $\bar{T}>0$ such that for every $T\geq\bar{T}$, the matrix
$\Phi_T$ is Hurwitz.
\begin{proof}
By assumption, $\lim_{T\rightarrow\infty}e^{AT}=0$. Therefore, and by (20),
\begin{equation}
\lim_{T\rightarrow\infty}\Phi_T=\left(\begin{array}{cc}
A & B\\
0 & -I_k
\end{array}
\right).
\end{equation}
By assumption, this matrix is Hurwitz.
\end{proof}
\end{lemma}

The next result concerns the tracking of a constant reference $r\in R^k$ by the output of an LTI system. Given $r\in R^k$. Fix $T>0$,
and assume that $A$ is nonsingular.
\begin{lemma}
Suppose that $\Phi_T$ is Hurwitz. Then
\begin{equation}
\lim_{t\rightarrow\infty}y(t)=r.
\end{equation}
\begin{proof}
By Eq. (22) with $r(t+T)\equiv r$, and the assumption that $\Phi_T$ is Hurwitz, the state variable $x(t)$ and the input $u(t)$
 have asymptotic values, namely,
there exist  $x\in R^n$ and $u\in R^k$ such that, $\lim_{t\rightarrow\infty}x(t)=x$ and $\lim_{t\rightarrow\infty} u(t)=u$.
Correspondingly, $\lim_{t\rightarrow\infty}\dot{x}(t)=0$ and $\lim_{t\rightarrow\infty}\dot{u}(t)=0$. Furthermore, by
the second part of Eq. (16), $\lim_{t\rightarrow\infty}y(t)=cx$.

Taking the limit $t\rightarrow\infty$, by (16), $Ax+Bu=0$. By (19) with $r(t+T)\equiv r$,
$r-Ce^{AT}x-CA^{-1}(e^{AT}-I_n)Bu=0$, and since $Ax+Bu=0$, we have that
\begin{equation}
r-Ce^{AT}x+CA^{-1}(e^{AT}-I_n)Ax=0.
\end{equation}
This implies, after some algebra, that $r-Cx=0$, hence Eq. (24) is satisfied.
\end{proof}
\end{lemma}

As a corollary of Lemma 1 and Lemma 2, if $A$ is Hurwitz then there exists $\bar{T}>0$ such that, with every $T\geq\bar{T}$, the system
will track any given constant $r\in R^k$.

\subsection{Example}
Consider the problem  of controlling the position of a particle by the force applied to it. The model we use is a second-order
system with a drag, obeying Newton's second law.
 Thus, denoting by
$y\in R$ the position of a particle with respect to a given reference, let $v\in R$ denote its velocity, and let $u\in R$ be the applied force.
The motion equations are
$\dot{y}(t)=v(t)$ and $\dot{v}(t)=av(t)+u(t)$, for given initial conditions $y_0:=y(0)$ and $v_0:=v(0)$, where $a\in R$ is the drag coefficient.
In a physical system $a<0$, but we allow for $a>0$ in order to extend the discussion to unstable
systems.  However, we assume henceforth that $a\neq 0$ since
the case where $a=0$ requires a different set of equations than those derived in the sequel. Given $r>0$, the objective is to regulate $y(t)$ to
$rt$, namely to achieve the limit $\lim_{t\rightarrow\infty}(y(t)-rt)=0$.

Defining the state variables by $x_1(t)=y(t)-rt$ and $x_2(t)=v(t)$, the state equation is
\begin{equation}
\begin{array}{cc}
\dot{x}_1(t)   = & x_2(t)-r\\
\dot{x}_2(t)  = & ax_2(t)+u(t),
\end{array}
\end{equation}
with initial conditions $x_{1,0}:=x_1(0)$ and $x_{2,0}:=x_2(0)$. The objective is to control $x_1(t)$ in the sense that
\begin{equation}
\lim_{t\rightarrow\infty}x_1(t)=0.
\end{equation}
By Eq. (13), and solving analytically the state equation (26), it is readily seen that
\begin{equation}
g(u(t))=x_1(t)+\frac{1}{a}(e^{aT}-1)x_2(t)+\dfrac{1}{a^2}(e^{aT}-1)u-\frac{T}{a}u-rT,
\end{equation}
and hence
\begin{equation}
\frac{\partial g}{\partial u}(u)=\frac{1}{a^2}(e^{aT}-1)-\frac{T}{a}.
\end{equation}
Considering Eqs. (15),  (26) with $r(t+T)\equiv 0$, and applying   Eqs. (14),  (28), (29), and (22) we obtain,   after some algebra, that
\begin{equation}
\Phi_{T}=\left(
\begin{array}{ccc}
0 & 1 & 0\\
0 & a & 1\\
-\frac{a^2}{e^{aT}-1-aT} & -\frac{a(e^{aT}-1)}{e^{aT}-1-aT} & -1
\end{array}
\right);
\end{equation}
the matrix $\Psi_{T}$ in (22) is irrelevant to the present discussion.
The characteristic polynomial of $\Phi_{T}$ is
\begin{eqnarray}
\chi(\lambda)=\lambda^3+(1-a)\lambda^2\nonumber \\
+\frac{a^2 T}{e^{aT}-1-aT}\lambda+\frac{a^2}{e^{aT}-1-aT}.
\end{eqnarray}
Lemma 2 ensures that   tracking, in the sense of Eq. (27),  is achieved if $\Phi_{T}$ is Hurwitz.

\begin{proposition}
The matrix $\Phi_{T}$ is Hurwitz if and only if the following two
conditions are satisfied: (i) $a<1$, and (ii) $T>\frac{1}{1-a}$.

\begin{proof}
Consider the Routh test.
Denote the entries of the first column in the Routh table by $c_{i}$, $i=3,2,1,0$, in decreasing  order.
Then $c_3=1$, $c_2=1-a$, $c_1=\frac{(1-a)a^2 T-a^2}{(e^{aT}-1-aT)(1-a)}$,
and $c_0=\frac{a^2}{e^{aT}-1-aT}$. Now $c_3>0$, while $c_2>0$ if and only if $a<1$. Next, for every $x\neq 0$, $e^{x}-1-x>0$,
this follows from the facts that, with the function $\zeta(x):=e^{x}-1-x$, $\zeta^{'}(x)>0$ for $x>0$,  $\zeta^{'}(x)<0$ for $x<0$, and $\zeta(0)=0$.
This implies that $c_0>0$. $\Phi_{T}$ is Hurwitz if and only if $a<1$ and $c_1>0$. For $a<1$, it is readily seen that $c_1>0$ if and only if $T>\frac{1}{1-a}$.
\end{proof}
\end{proposition}

Note that  for $a>1$ no $T>0$ will yield tracking.   On the other hand,
for $a\in(0,1)$,  tracking is attained as long as $T>\frac{1}{1-a}$ even though the plant is unstable.

We mention that for $a>1$ it is possible to scale up the RHS of Eq. (15)  in order to guarantee tracking.
Multiplying the RHS of Eq. (15) by $\alpha>1$  results in the scaling of the last row of the matrix
$\Phi_{T}$ in  (30)
by $\alpha$. A bit of algebra reveals that the entries in the first column of the Routh table are
$c_3=1$, $c_2=\alpha-a$, $c_1=\frac{(\alpha-a)\alpha a^2T-a^2 \alpha}{(e^{aT}-1-aT)(\alpha-a)}$,
and $c_0=\frac{\alpha a^2}{e^{aT}-1-aT}$. The arguments in the proof of Proposition 2 yield that tracking is attained if and only if $a<\alpha$ and $T>\frac{1}{\alpha-a}$.
All of this will be demonstrated by simulation in the next section.

\section{Simulation and Laboratory Experiments}
This section presents simulation results on three systems: the position-control system described in the last section,
an inverted pendulum, and a  platoon of mobile robots. For the first system we consider the tracking of a ramp
where we verify the theoretical results derived in Section II.D, then we consider an
unstable system tracking a sinusoid. The inverted pendulum provides an example of a nonlinear, unstable system  tracking a sinusoid,
and in the platoon example we control the interspacing between successive vehicles.
The first two systems are simulated by MATLAB codes for solving their respective state equations via Euler's
forward method with integration step sizes of $dt=0.01$. In contrast, the platoon system is implemented in a laboratory setting. In all three cases the
lookahead simulation for computing $g(u)$ and $\frac{\partial g}{\partial u}(u)$, hence the control law, involves the Forward Euler method  with integration
step sizes of $\Delta t:=0.01T$.

\subsection{Position-control system}
Consider the system described in Section II.D, where the objective is to have the position $y(t)$ track a ramp $rt$ for
$r=2$. The system is defined by Eq. (26), and the objective is to attain Eq. (27). The controller is defined by (15), and by
Lemma 2, tracking is attained as long as the matrix $\Phi_{T}$ is Hurwitz. By Proposition 2, this is the case if and only if
$a<1$ and $T>\frac{1}{1-a}$. We next verify this conclusion.

Let $a=-1$ so that the plant system is stable, and let  $T=1$ which ensures that the closed-loop system is stable.
Following a simulation of the system,  the graph of $x_1(t)$ vs. $t$ is depicted in Figure 2 by the solid graph, and it indicates tracking.
Next we change $T$ to $T=0.5$, which is on the boundary of the stability region. The resulting graph of $x_1(t)$,  shown by
the dashed graph in Figure 2, indicates growing  oscillations and hence instability. Theoretically one may expect undamped oscillations
at $T=0.5$
because the closed-loop system is marginally stable, but that is obtained for a slightly larger $T$, namely $T=0.53$ (results not shown).
The discrepancy is due to
the use of the forward Euler integration which can destabilize marginally-stable systems. Several experiments with various values of $T<0.53$,
not shown here, resulted in instability of the closed-loop system.

For the next experiment we make two changes to the system: $a$ is set to  $a=0.5$ hence the plant system is unstable, and the tracking's target
is a time-varying signal $r(t):=2+sin~t$. Due to the value of $a$, stability of the closed-loop system requires that $T>2$.
First we take $T=3$ to ensure stability.  The resulting graph of the position $y(t)$  is depicted by
the dashed curve in Figure 3, and for the sake of reference we plotted $r(t)$ by the  dot-dashed curve.\footnote{Actually the dot-dashed curve  looks like a solid
curve after the first half-cycle of the sinusoid since it is merged there with another graph (described below). The dot-dashed part of the sinusoid
is visible during its first half-cycle.} We notice stability of the closed-loop system but no tracking; larger
simulation horizons do not change this conclusion. It  is not surprising in light of the discussion in Section II,
since $T$ is too large for $g(u)$ to yield an adequate approximation for $y(t+T)$. Therefore we reduce
$T$ to $T=0.4$ in order to get a better approximation. However, this value is not in the stability range. Therefore we speed up the
 controller by multiplying the RHS of Eq. (15) by $\alpha=5.0$.
 The resulting graph of $y(t)$ is depicted by the solid curve in Figure 3, and we see there
stability and as well as tracking. It must be pointed out that in the first 3 seconds the motion may experience
a large jerk. This issue can be addressed by ad-hoc methods like gradual increases of the gain $\alpha$, which are beyond the
scope of this paper but will be considered in the near future.

\vspace{.2in}
\begin{figure}[h]
\centering
\vspace{-.3cm}
\includegraphics[width=0.40\textwidth]{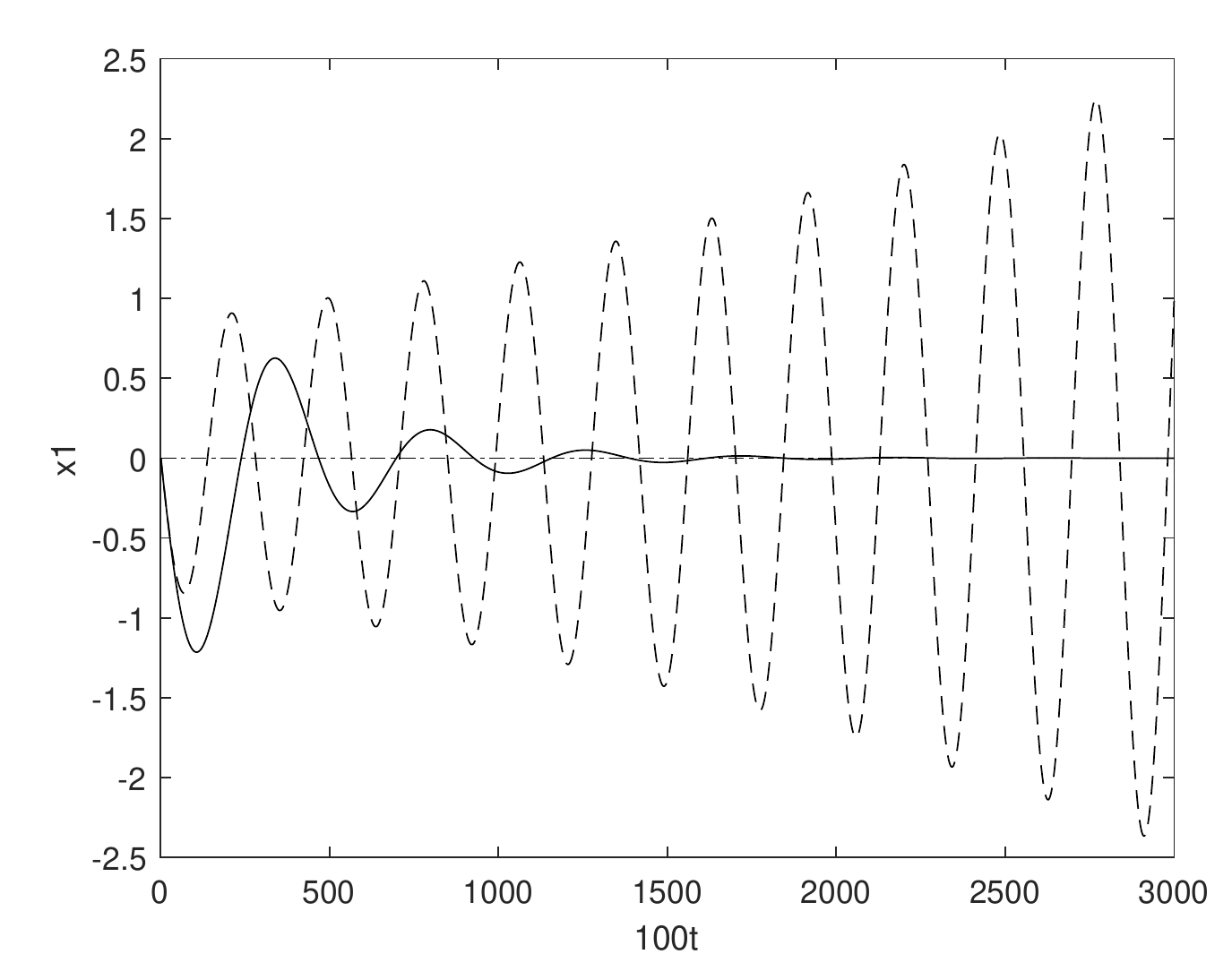}
% %\includegraphics[width=3.15in,angle=0]{Adaptive.jpg}
{\small \caption{Position-control system: constant target}}
\end{figure}

\vspace{.2in}
\begin{figure}[h]
\centering
\vspace{-.3cm}
\includegraphics[width=0.40\textwidth]{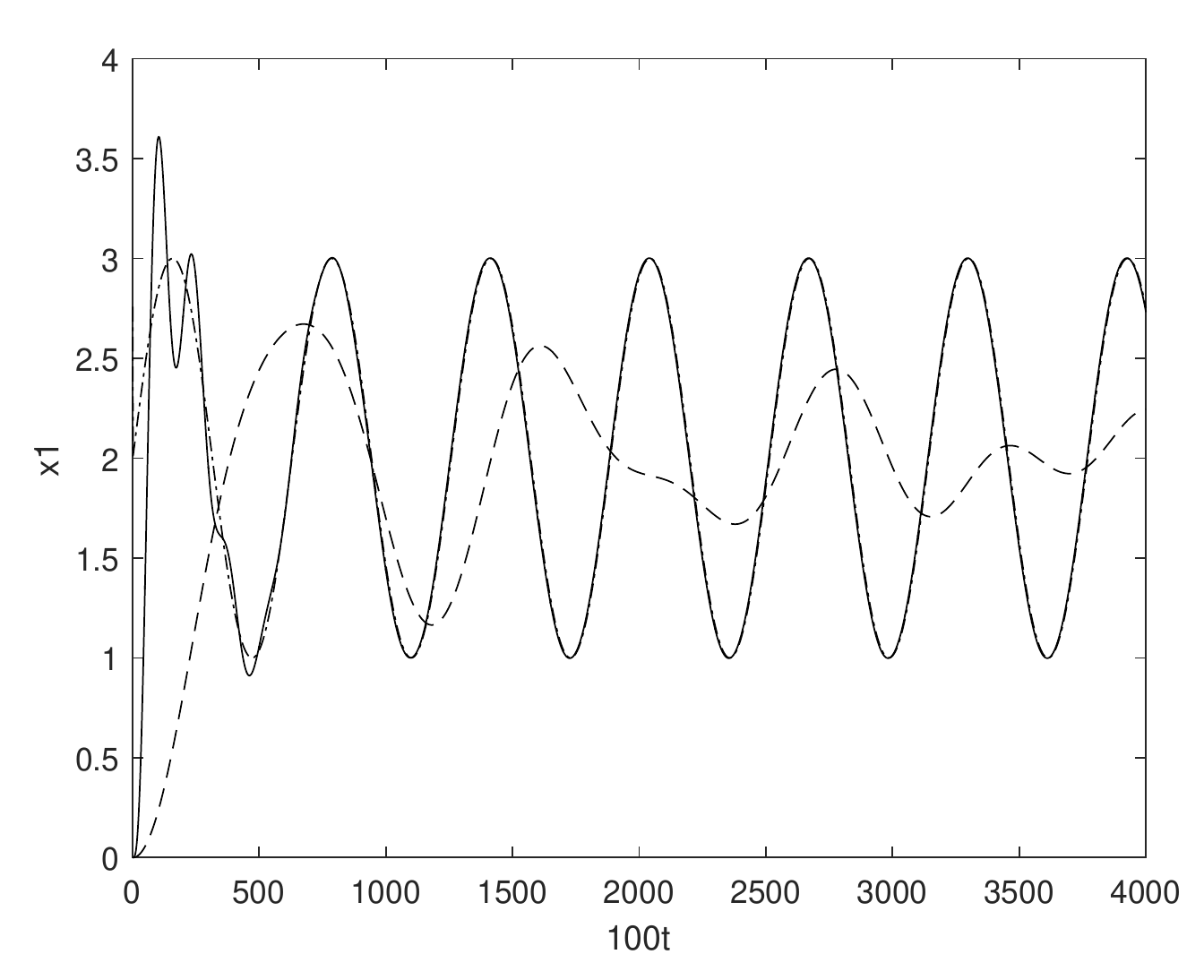}
% %\includegraphics[width=3.15in,angle=0]{Adaptive.jpg}
{\small \caption{Position control system: sinusoidal target}}
\end{figure}

\subsection{Inverted pendulum}

The equations for the inverted-pendulum that we use are
\begin{equation}
\begin{array}{cc}
\dot{x}_1(t)   = & x_2(t)\\
\dot{x}_2(t)  = & a~sin~x_1~-bx_2 +u(t),
\end{array}
\end{equation}
where $x_1$ is the angle of the pendulum from the upper equilibrium point and  $x_2$ is its angular velocity; $a>0$ and $b>0$ are given positive
constants. We chose (arbitrarily) $a=1.0$ and $b=0.2$.
In the first experiment we regulate the angle $x_1$ to the target $\pi/6$. We set the lookahead parameter $T$ first to $T=2.0$ and then to $T=0.8$.
The results are
shown in Figure 4, where the dotted horizontal line indicates the  target level of  $\pi/6$.
For $T=2$, the graph of the angle $x_1(t)$ is depicted by the solid curve, and we discern stability and tracking. However, for the smaller
value of $T=0.8$, the graph of $x_1$, depicted by the dashed curve,  indicates that the closed-loop system is unstable and  no tracking is achieved. These results are consistent with the discussion in Section
II suggesting that stability of the closed-loop system is more likely to be attained for larger values of $T$.

\vspace{.2in}
\begin{figure}[h]
\centering
\vspace{-.3cm}
\includegraphics[width=0.40\textwidth]{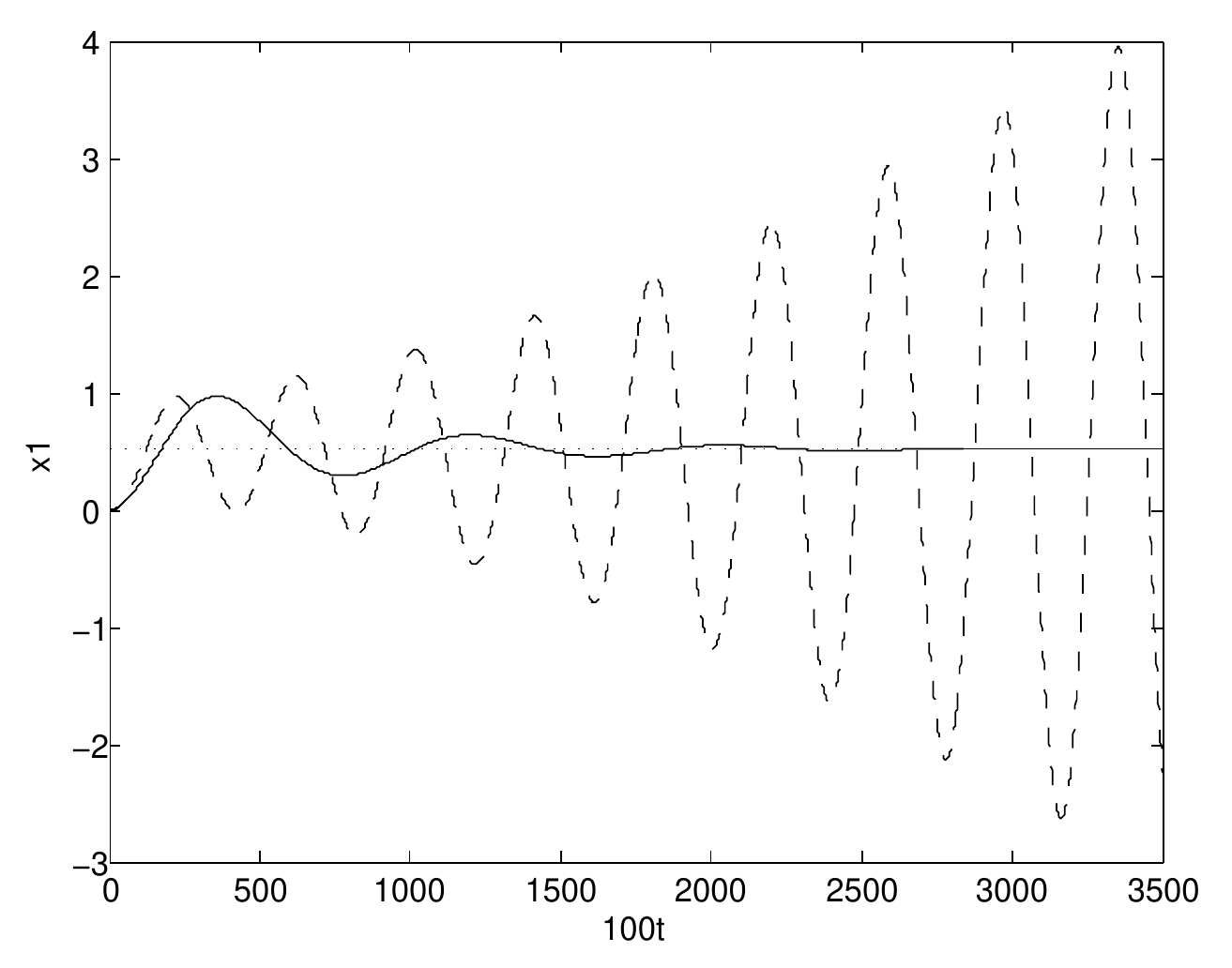}
% %\includegraphics[width=3.15in,angle=0]{Adaptive.jpg}
{\small \caption{Inverted pendulum: constant target}}
\end{figure}

Consider next the tracking of the curve $r(t)=\frac{\pi}{6}+\frac{\pi}{8}sin~t$. For $T=2.0$, the results are shown in Figure
5, where the dash-dotted curve  depicts the graph of $r(t)$, and the dashed graph is of $x_1(t)$.\footnote{It may
 be hard to distinguish the dash-dotted curve from a solid curve with which
 it is almost aligned. The solid curve will be explained shortly.} It is evident that the closed-loop system
is stable but no tracking is achieved.
 We then set $T$ to $T=0.15$, but the resulting closed-loop system is unstable; in fact, the oscillations in
$x_1(t)$ (not shown in the figure) reach a magnitude of about $10^{21}$ at $t=35$.  To address this problem
we speed up the controller by multiplying the RHS of (15) by  $\alpha=20$. The resulting graph of $x_1(t)$ is depicted in Figure 5 by the solid curve, and
it is evident that tracking has been attained. We point out that the parameter-values $T=0.15$ and $\alpha=20$ were chosen after
some trial and error, but during that process it was quite evident that smaller $T$ yields better tracking but requires
larger $\alpha$. For instance, increasing $T$ from $T=0.15$ to $T=0.2$, $\alpha=8$ sufficed to give tracking and the graph of
$x_1(t)$ was quite similar to that shown in Figure 5. Judging by the error ${\cal E}:=\int_{5}^{35}|r(t)-x_1(t)|dt$,
for $T=0.15$ and $\alpha=20$, ${\cal E}=1.056$; while for $T=0.2$ and $\alpha=8$, ${\cal E}=1.419$.
These errors correspond to  average errors $|x_1(t)-r(t)|$ over $t\in[5,35]$ of $0.035$ and $0.047$, respectively.

\vspace{.2in}
\begin{figure}[h]
\centering
\vspace{-.3cm}
\includegraphics[width=0.40\textwidth]{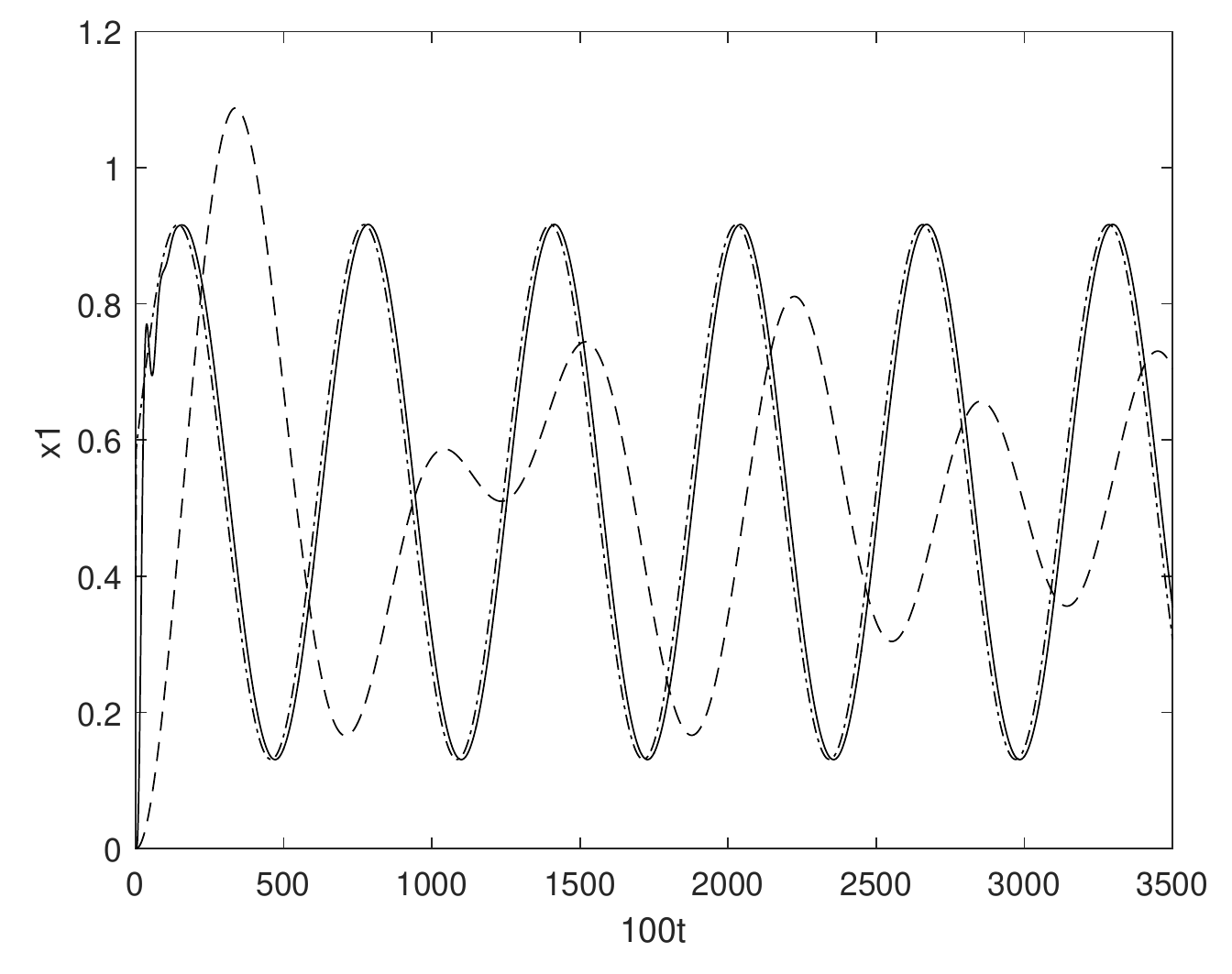}
% %\includegraphics[width=3.15in,angle=0]{Adaptive.jpg}
{\small \caption{Inverted pendulum: sinusoidal  target}}
\end{figure}

\subsection{Platoon of mobile robots}

Consider a platoon of  mobile robots (vehicles) tasked with sequentially following a given path in $R^2$ in a predetermined order.
Denote by $R_{i}$, $i=1,\ldots,N$ the robots in their sequential order.
The vehicles are autonomous in the sense that each one of them controls its own motion. The objective is to have the first (leading) robot, $R_1$,
regulate its speed to a given reference, and for every $i=2,\ldots,N$, $R_i$, is to regulate its distance to $R_{i-1}$.
This problem has been extensively investigated in the context of automated highway and urban traffic control,  and recently an interest in it
has sprouted in the setting of smart cities; see. e.g., \cite{Hunter13, Cassandras16}  for  surveys.

The experiments described below were run in a  swarm robotic testbed laboratory,  the
{\it Robotarium}, situated in the Georgia Tech campus\cite{Pickem17}. The robots, GRITSbots \cite{Pickem15}, and their motion
can be controlled directly by their velocities via differential-drive motors.
The differential drive robots are modelled by unicycle dynamics. Following \cite{Olfati_Saber02}, rather than directly specifying the translational and angular velocities of each robot, we instead control a point-particle in the front of  of each robot with
the simple  dynamics,
\begin{equation}
\dot{x}_{i}=u_i,
\end{equation}
and map this onto the dynamics of the robot.

The subject of the experiment is a platoon of    eight robots, and  the goal is to have them move counter-clockwise on a given reference circle at a
 predetermined speed while maintaining a given  interspacing of $d$ cm between $R_{i}$ and $R_{i-1}$.
The center and radius of the reference circle are
given  $C\in R^2$ and $r>0$, hence it is denoted by $B(C,r)$.   The first robot, $R_1$, starts on $B(C,r)$
and is programmed to stay on it while controlling   its velocity to the given target.   Each subsequent robot $R_i$, $i=2,\ldots,N$, attempts, at time $t$,  to regulate its  position  to the point of (Euclidean) distance $d$ behind $R_{i-1}$ on the circle $B(C,r)$.
The parameters for the experiments are as follows: $r=28$cm and   $d=14$cm, the robot's dimension is $3cm\times 3cm$, and the point  which is controlled
is $3cm$ ahead of the robot. The  lookahead parameter for the controller (Eq. (15)) is
T=0.6s, and the scaling factor of the RHS of (15) is $\alpha=45$.

The results of the experiments can be seen in the video clip \cite{Buckley17}, where we see that although all the robots start close to  the circle $B(C,r)$, some of their trajectories initially move away from it  and display an erratic behaviour. However, they soon turn to the circle
 and track their target distance from each other.
 Figure~6 displays the progress of the tracking assignments by providing  snapshots of the robots' locations at times $t=0$,  $t=7$s (200 iterations), and $t=20$s (600 iterations). Finally,
Figure~7 depicts the graphs of the distance between adjacent robots,  $||x_{i}-x_{i-1}||$,  with roughly 30
iterations per second, and we see that the interspace  tracking has been achieved.

\vspace{.2in}
\begin{figure}[h]
\centering
\vspace{-.3cm}
\includegraphics[width=0.4\textwidth]{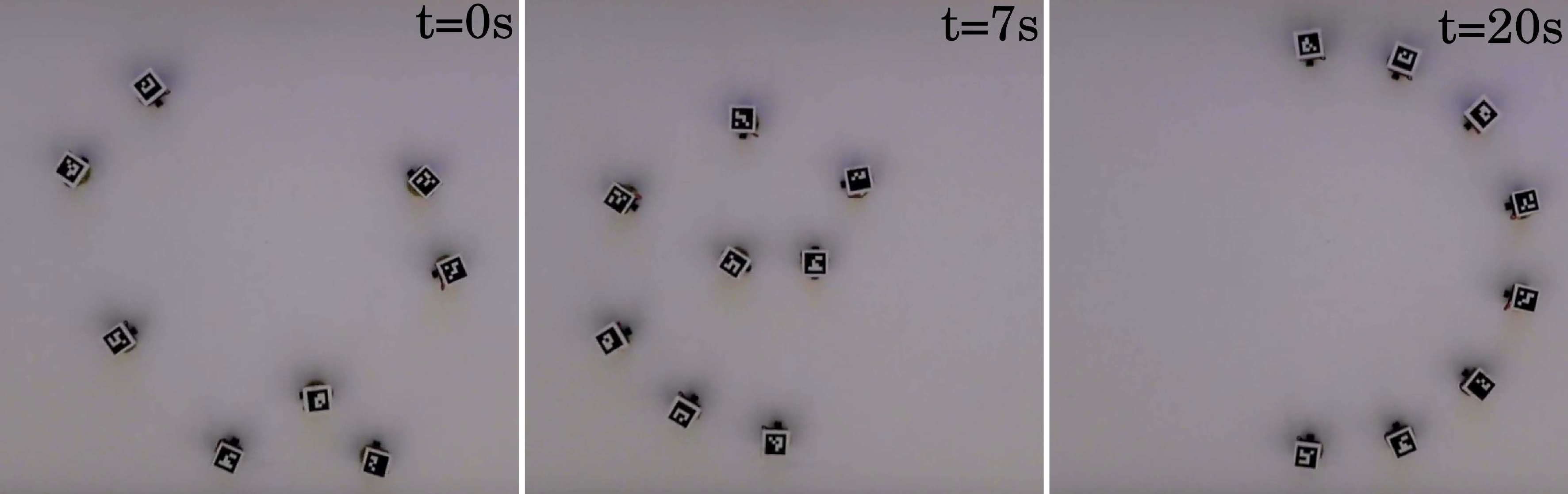}
% %\includegraphics[width=3.15in,angle=0]{Adaptive.jpg}
{\small \caption{Positions of robots at times at times 0 (left), 7s (center),
and 20s (right)}}
\end{figure}

\vspace{.2in}
\begin{figure}[h]
\centering
\vspace{-.3cm}
\includegraphics[width=0.40\textwidth]{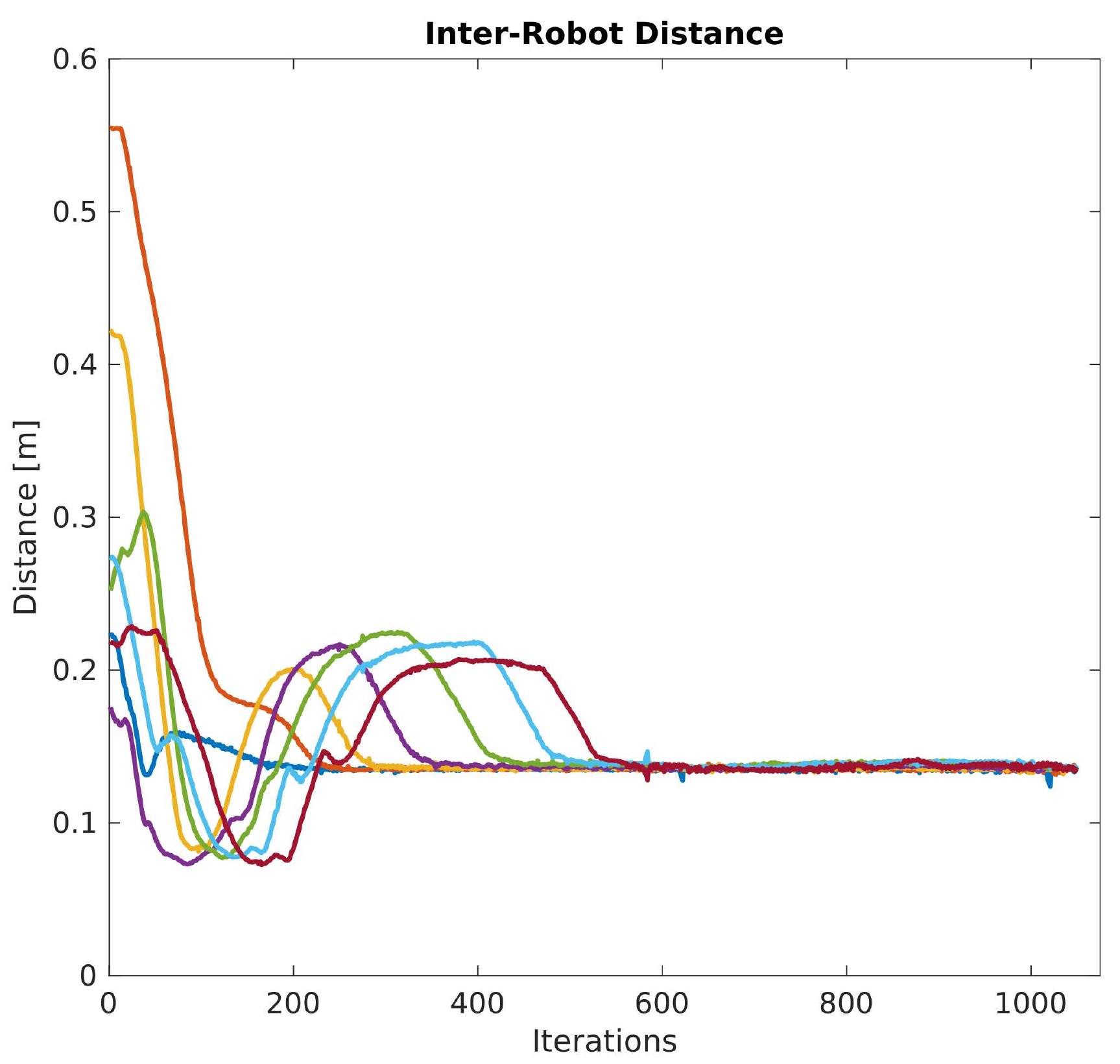}
% %\includegraphics[width=3.15in,angle=0]{Adaptive.jpg}
{\small \caption{Interspacing between adjacent robots}}
\end{figure}

\section{Conclusions}

This paper proposes a technique for performance regulation and tracking of dynamical systems via lookahead simulation.
It is based on steering the  control signal in a direction determined by the Newton-Raphson method for solving an algebraic
equation, the loop equation. The resulting controller can be nonlinear. Preliminary theoretical results are derived for
memoryless nonlinear systems and dynamic linear systems. Simulation experiments support the theoretical developments, and go
beyond them to include nonlinear  dynamical systems such as the inverted pendulum. A laboratory experiment of interspacing control
in a platoon of mobile robots is described as well.

A key question concerns the determination of the lookahead timing
parameter $T$ for the simulations defining the controller. The derived theoretical and simulation results indicate that stability of the closed-loop system may require large $T$, while tracking may require small $T$. This tradeoff  has been resolved for
the systems under study by first picking a small $T$, then stabilizing the system (if needed) by
speeding up the action of the controller. Identifying the kind of systems for which this regulation and tracking technique works constitutes a subject of future research.

%\bibliographystyle{IEEEtran}
%\bibliography{biblio}

\end{document}